\documentclass{amsart}

\usepackage[latin1]{inputenc}

\usepackage{amsmath,amssymb,amsthm}

\usepackage[english]{babel}

\def\FPt{{\sf FPt}}
\def\WPt{{\sf WPt}}
\def\diag{{\sf diag}}
\def\sign{{\sf sign}}
\newcommand{\Sld}{\ensuremath{\mathcal{S}(\mathbb{L}^d})}
\newcommand{\Int}{\text{Int}}
\newcommand{\Ext}{\text{Ext}}
\newcommand{\Color}{\text{Color}_{\phi,\nu}^r}
\newcommand{\Intbarre}{\overline{\text{Int}}}
\newcommand{\Extbarre}{\overline{\text{Ext}}}
\newcommand{\peierls}{\mathcal{P}\text{eierls}}
\newcommand{\as}{a.s.}

\newcommand{\N}{\mathbb{Z}_{+}}

\newcommand{\Z}{\mathbb{Z}}
\newcommand{\Zp}{\mathbb{Z}_{*}^2}
\newcommand{\Zd}{\mathbb{Z}^d}

\newcommand{\R}{\mathbb{R}}
\newcommand{\Rd}{\mathbb{R}^d}

\newcommand{\E}{\mathbb{E}\ }
\newcommand{\EE}{\mathbb{E}}
\newcommand{\Ld}{\mathbb{L}^d}
\newcommand{\Ldeux}{\mathbb{L}^2}

\newcommand{\Ldeuxstar}{\mathbb{L}_{*}^2}
\newcommand{\Ed}{\mathbb{E}^d}
\newcommand{\Fd}{\mathbb{F}^d}
\newcommand{\Edeux}{\mathbb{E}^2}
\newcommand{\Edeuxstar}{\mathbb{E}^2_{*}}

\newcommand{\Cov}{\text{Cov}}
\newcommand{\ie}{\emph{i.e. }}
\newcommand{\eg}{\emph{e.g. }}
\newcommand{\resp}{\emph{resp. }}

\newcommand{\miniop}[3]{%
\renewcommand{\arraystretch}{0.6}
\begin{array}{c}
{\scriptstyle #1}\\
#2\\
{\scriptstyle #3}
\end{array}
\renewcommand{\arraystretch}{1}}

\newcommand{\Card}[1]{\vert #1 \vert}

\newcommand{\communique}{\leftrightarrow}

\newcommand{\1}{1\hspace{-2.7mm}1}

\newtheorem{theorem}{Theorem}
\newtheorem{lemme}{Lemma}

\renewcommand{\P}{\mathbb{P}}
\newcommand{\bor}[1][]{\mathcal{B}(#1)}

\newtheorem{prop}{Proposition}

\begin{document}

\title{Central limit theorems in Random cluster and Potts Models}
\author{Olivier Garet}
\date{\today}
\address{Laboratoire de Mathématiques, Applications et Physique
Mathématique d'Orléans UMR 6628, Université d'Orléans, B.P.
6759,
 45067 Orléans Cedex 2 France}
\email{Olivier.Garet@labomath.univ-orleans.fr}

\subjclass{60K35, 82B20, 82B43.} 

\keywords{random cluster measure, percolation, coloring model, Central Limit Theorem, Potts Model, Ising Model.}

\begin{abstract}
We prove that for $q\ge 1$, there exists $r(q)<1$ such that for $p>r(q)$, the number of points in large boxes which belongs to the infinite cluster has a normal central limit behaviour under the random cluster measure $\phi_{p,q}$ on $\Zd$, $d\ge 2$. Particularly, we can take $r(q)=p_g^{*}$ for $d=2$, which is commonly conjectured to be equal to $p_c$. These results are used to prove a $q$-dimensional central limit theorems relative to the fluctuation of the empirical measures  for the ground Gibbs measures of the $q$-state Potts model at very low temperature and the Gibbs measures which reside in the convex hull of them. A similar central limit theorem is also given in the high temperature regime. Some particular properties of the Ising model are also discussed.  
\end{abstract}

\maketitle

\section{Introduction}
The aim of this study is to answer to a natural question relative to Gibbs measures in the $q$-state Potts model: take a finite box in $\Zd$ and consider the frequencies of occurring of each of the $q$ state in the box. It is obvious that the vector of empirical frequencies converges to a constant when the considered Gibbs measure is ergodic. Now,  several natural naturally arise. Two years ago, Cerf and Pisztora considered the difficult problem of large deviations~\cite{MR2003g:82018}.
We will  consider here the problem of having a central limit theorem with a standard renormalization.
To this aim, we will use the road of the  Fortuin-Kasteleyn random cluster measure. In the last decade, consequent progresses have been made in the study of this model -- see Grimmett~\cite{Grimmett-review} for a large panorama --, and  it also appeared that in most cases, the shortest way to results relative to the ferromagnetic Potts model encounters random cluster model -- see the review of H{\"a}ggstr{\"o}m~\cite{MR2000e:60182} for a self-contained introduction to these relations. 
 
Roughly speaking, we can say that a realization of the $q$-state Potts model
with free boundary conditions in a finite box is a random coloring of the vertices of a realization of
a free random-cluster measure in the box $\Lambda$, with the constraint that
connected components are mono-color.
 Actually, we can consider the Potts model as the restriction to its vertices of a measure on ``colored graphs'': there is randomness  on the set of open bonds and also on the color of vertices, with the condition that connected components are mono-color.
Biskup, Borgs, Chayes and Koteck{\'y} have shown that this approach could also be fruitful in the study of Potts models with external fields~\cite{MR2001c:82015}. It leads them to introduce a generalization of the FK random cluster measure. 

It is not difficult to guess that the presence (or the absence) of an infinite cluster strongly modifies the fluctuation of the empiric repartitions.
So it appears that the study of the supercritical case -- which is, of course, the most interesting one -- necessitates some knowledge about the normal fluctuations of the size of the intersection of the infinite cluster with large boxes. We made some progress in this aim but did not success to have such a result in the whole supercritical region. This gap in the spectrum of results, and
the intuition we have that the random coloring scheme had already much to bring, lead us to adopt the following organization for our paper.
\begin{itemize}
\item In a first part, we prove that for each   $q\ge 1$, there exists $r(q)<1$ such that for $p>r(q)$, the number of points in large boxes which belongs to the infinite cluster has a normal central limit behaviour under the random cluster measure $\phi_{p,q}$. The result is much better on the square lattice: in
this case, we can take $r(q)=p_g^{*}$, which is commonly conjectured to be equal to $p_c$ -- it is even proved for $q=1,2$ and $q\ge 22$.
\item In a second part, we prove a $d$-dimensional central limit theorem
for the fluctuation of the empiric repartitions of colors in a coloring model, that is model where the connected components of a random graph are paint
independently. The random graph measure is supposed to be ergodic, to have appropriate moments for the size of finite clusters and, in case it has infinite cluster, to satisfy to a central limit theorem for the fluctuations of the size of the intersection of the infinite cluster with large boxes. This generalizes and extends the result obtained by the author~\cite{MR2003g:60166} for the fluctuations of the magnetization in the case where the random graph measure is the classical Bernoulli measure.
\item In a third part, we combined these results to obtain central limit theorems for the $q$-state Potts model which may be summarized as follows: When $\beta<\beta_g=-\frac12\ln(1-p_g)$, the vector of empiric repartitions satisfy to  an ordinary central limit theorem under the unique Gibbs measure.
If  $\beta>\beta_r$, we have a central limit theorem with a random centering under each Gibbs measure which is in the convex hull of the ground Gibbs measures. The limit is normal for the ground Gibbs measures, but not necessarily for their mixtures. The result is also better for $d=2$ because $\beta_r=\beta_c$ as soon as $r(q)=p_c$.
In the case of the Ising model, the limit is always Gaussian even for a mixture
of ground Gibbs measures, so we can say that for $d=2$ and $\beta\ne\beta_c$,
the fluctuation of the empirical magnetization around their (random)
limit magnetization are asymptotically normal, whatever Gibbs measure we consider.

\end{itemize}     

\section{Notations and preliminaries}

\noindent\textbf{Graph theoretical notations}\\ 
For $x\in\Zd$, let us denote $\Vert x\Vert=\sum_{i=1}^d | x_i|$ and
now consider the 
graph $\Ld=(\Zd,\Ed)$, with
$$\Ed=\{\{x,y\}; x,y\in\Zd\text{ and }\Vert x-y\Vert=1\}.$$

For $x\in\Zd$ and $r\in [0,+\infty)$, we note $B(x,r)=\{y\in\Zd;\Vert x-y\Vert\le r\}$. 
If $e=\{x,y\}\in\Ed$, then $x$ and $y$ are called \textit{neighbours}.

 In the following, the
expression "subgraph of $\Ld$" will always be employed for each
graph of the form $G=(\Zd,E)$ where $E$ is a subset of $\Ed$. We denote by 
$\Sld$ the set of all subgraphs of $\Ld$.

Set $\Omega=\{0,1\}^{\Ed}$. An edge $e \in \Ed$ is said to be
\textit{open} in the configuration $\omega$ if $\omega(e)=1$, and \textit{closed}
otherwise. 

There is a natural bijection between $\Sld$ and $\Omega$, that is $E\mapsto (1_{e\in E})_{e\in\Ed}$.
Consequently, we sometimes identify $\Sld$ and $\Omega$ and say ``random graph measure'' rather than ``measure on $\Omega$''.

A \textit{path} is a sequence $\gamma=(x_1,
e_1,x_2,e_2,\ldots,x_n,e_n,x_{n+1})$ such that $x_i$ and $x_{i+1}$ are
neighbours and $e_i$ is the edge between $x_i$ and $x_{i+1}$. 
We will also sometimes describe $\gamma$ only by the vertices it
visits $\gamma=(x_1,x_2,\ldots,x_n,x_{n+1})$ or by its edges
$\gamma=(e_1,e_2,\ldots,e_n)$. The number $n$ of edges in $\gamma$ is called
the \textit{length} of $\gamma$ and is denoted by $|\gamma|$. Of particularly
interest are the \textit{simple paths} for which the visited vertices are all
distinct. 
We will also consider \textit{cycles}, that are paths for which the visited vertices are all distinct, except that $x_1=x_{n+1}$. 
A path is said to be \textit{open} in the configuration $\omega$
if all its edges are open in $\omega$.

The \textit{clusters} of a configuration $\omega$ are the connected
components of the graph induced on $\Z^d$ by the open edges in
$\omega$. For $x$ in $\Z^d$, we denote by $C(x)$ the cluster containing
$x$. In other words, $C(x)$ is the set of points in $\Z^d$ that are linked
to $x$ by an open path. 
We note $x\communique y$ to signify that $x$ and $y$ belong to the same cluster. If $A$ and $B$ are subsets of $\Zd$, $A\communique B$ means that there exists
$(x,y)\in A\times B$ with $x\communique y$.
We also note $x\communique\infty$ to say that $|C(x)|=+\infty$.
In the whole paper, we will note $I=\{x\in\Zd; x\communique\infty\}$. 

We say that two bonds $e$ and $e'$ of $\Ed$ are neighbours if
$e\cap e'$ is not empty. It also gives a notion of connectedness in $\Ed$ by
the usual way.
 
For each subset $\Lambda$ of $\Zd$, we will note by $\partial\Lambda$ the boundary of $\Lambda$:
$$\partial\Lambda=\{y\in\Lambda^c;\exists x\in\Lambda\text{ with }\Vert x-y\Vert=1\}$$
and $\EE_{\Lambda}$ the set of inner bonds of $\Lambda$:
$$\EE_{\Lambda}=\{e\in\Ed; e\subset \Lambda\}.$$ 
Note that if $\Lambda$ and $\Lambda'$ are disjoint sets, then
$\EE_{\Lambda}$ and $\EE_{\Lambda'}$ are disjoint too.

For each $E\subset \Ed$, we denote by $\sigma(E)$
the $\sigma$-field generated by the projections $(\omega_e)_{e\in E}$.\\ 
When $\Lambda\subset\Zd$, we also use the notation $\sigma({\Lambda})$
instead of $\sigma(\EE_{\Lambda})$.\\

We sometimes consider another set of bonds on $\Zd$, that is
$$\Fd=\{\{x,y\}; x,y\in\Zd\text{ and }\Vert x-y\Vert_{\infty}=1\},$$
where $\Vert x\Vert_{\infty}=\max(|x_i|;1\le i\le d)$.
If $e=\{x,y\}\in\Ed$, then $x$ and $y$ are called $*$-\textit{neighbours}.
Similarly, we define the notion of $*$-paths, $*$-cycles, $*$ connected sets,\dots exactly in the same way that for the graph $\Ld$. 

For our central limit theorems, we will use boxes $(\Lambda_t)_{t\ge 1}$, with
$$\Lambda_t=\{x\in\Zd;\Vert x\Vert_{\infty}\le t\}.$$

Let $X$ and $S$ be arbitrary sets. 
Each ${\omega}\in X^{S}$ can be considered as a map from $S$ to $X$.
We will denote ${\omega}_{\Lambda}$ its restriction to $\Lambda$.
Then, when $A$ and $B$ are two disjoint subsets of $S$ and $({\omega},\eta)\in X^A\times X^B$,
${\omega}\eta$ denotes the concatenation of ${\omega}$ and ${\eta}$, that is the element $z\in X^{A\cup B}$ such that
\begin{equation*}
z_{i}=
\begin{cases}
\omega_i & \text{if }i\in A\\
\eta_i & \text{if }i\in B.
\end{cases}
\end{equation*}

\noindent\textbf{FKG inequalities}\\ 
If $\phi$ is a probability measure and $f,g$ two measurable functions, we note
$$\Cov_{\phi}(f,g)=\int fg\ d\mu-\big(\int fg\ d\mu\big)\big( \int fg\ d\mu\big).$$ If $A$ and $B$ are measurable events, be also note $\Cov_{\phi}(A,B)=\Cov_{\phi}(\1_A,\1_B)=\phi(A\cap B)-\phi(A)\phi(B)$.

If $X$ in an ordered set, we say that a measure $\phi$ on  $X^{S},\bor[X^S]$ satisfy to the FKG inequalities if for each pair of increasing functions $f$ and $g$, we have $\Cov_{\phi}(f,g)\ge 0$.

An event $A$ is said to be increasing (\resp decreasing) if $\1_A$ (\resp $1-\1_A$) is an increasing function. Of course, if $\phi$ satisfy to the  FKG inequalities and $A$ and $B$ are increasing events, we have $\Cov_{\phi}(A,B)\ge 0$.

\subsection{FK Random cluster measures}
Let $0\le p\le 1$ and $q>0$.

For each configuration  $\eta\in\Omega$ and each connected subset $E$ of $\Ed$
 we define the random-cluster measure $\phi_{E,p,q}^\eta$ with boundary condition $\eta$ on $(\Omega,\mathcal{B}(\Omega))$  by 
$$
\phi_{E,p,q}^\eta (\omega ) = 
\begin{cases}\dfrac 1{Z_{\Lambda,p,q}^\eta}
\bigg\{\displaystyle\prod_{e\in E} p^{\omega (e)} (1-p)^{1-\omega (e)}\bigg\}
q^{k_{\Lambda}(\omega_{\Lambda}\eta_{\Lambda^c} )} &\text{if $\omega_{E^c}=\eta_{E^c}$},\\
0 &\text{otherwise},
\end{cases}
$$
where $k(\omega)$ is the number of components of the graph
in the configuration $\omega$ which intersect $E$.
$Z_{E,p,q}^\eta$ is the renormalizing constant
$$
Z_{E,p,q}^\eta =\sum_{\omega\in\{0,1\}^{E}}\bigg\{\prod_{e\in E}
p^{\omega (e)} (1-p)^{1-\omega (e)}\bigg\} q^{k(\omega_{E}\eta_{E^c})}.
$$

For each $b\in\{0,1\}$, we will simply denote by $\phi_{E,p,q}^b$ the measure $\phi_{E,p,q}^\eta$ corresponding to the configuration $\eta$ which is such that $\eta_e=b$ for each $e\in\Ed$.

When $\Lambda$ is a finite subset of $\Zd$, we also use the notation $\phi_{\Lambda,p,q}$ instead of  $\phi_{\EE_\Lambda,p,q}$

A probability measure $\phi$ on $(\Omega,\bor[\Omega])$ is called a {\it random-cluster measure\/} with parameters $p$ and
$q$ if for each measurable set $A$ and each finite subset
$\Lambda$ of $\Zd$, we have the D.L.R. condition:
$$\phi(A)=\int_{\Omega}\phi_{\Lambda,p,q}^{\eta}(A)\ d\phi(\eta)$$
or equivalently if for each finite subset
$\Lambda$ of $\Zd$, we have the equation
$$\phi (.|\sigma(\Ed\backslash\EE_{\Lambda}))(\eta)=\phi_{\Lambda,p,q}^\eta\quad\phi -\as$$
The set of such measures is denoted by $\mathcal{R}_{p,q}$.

Let $b\in\{0, 1\}$. If $(\Lambda_n)_{n\ge 1}$ is an increasing sequence of volumes tending to $\Zd$,
it is known that the sequence $\phi_{\Lambda_n,p,q}^b$ as a weak limit which does not depend of the sequence $(\Lambda_n)_{n\ge 1}$. We denote by $\phi_{p,q}^b$ this limit. The following facts are well known; refer to the recent summary of Grimmett~\cite{Grimmett-review} of  complete references.
\begin{itemize}
\item $\phi_{p,q}^b$ is a translation  invariant ergodic measure.
\item $\phi_{p,q}^b\in \mathcal{R}_{p,q}$.
\item If $q\ge 1$, $\phi_{p,q}^b$ satisfies to the FKG inequalities.
\item Let us note $\theta^b(p,q)=\phi_{p,q}^b(0\communique\infty)$. 
There exists $p_c(q)\in (0,1)$, such that for each $b\in\{0,1\}$ we have  $\theta^b(p,q)=0$ for $p<p_c(q)$ and $\theta^b(p,q)>0$ for $p>p_c(q)$.
\end{itemize}

\noindent\textbf{Exponential bonds}\\ 
Grimmett and Piza~\cite{MR99a:60115} also introduced another critical probability:
Let us define
$$
Y(p,q)=\limsup_{n\to\infty}
\Bigl\{ n^{d-1}\phi_{p,q}\bigl( 0\communique\partial B(0,n)\bigr)\Bigr\}
$$
and
$p_g(q)=\sup\bigl\{ p:Y(p,q)<\infty\bigr\}$.
We have $0<p_g(q)\leq p_c(q)$, and it is believed that 
$p_g(q)=p_c(q)$ for all $q\ge 1$. 
They proved the following exponential bound:

\begin{prop}
\label{gr-pi}
Let $q\ge 1$, $d\ge 2$. For $p<p_g(q)$, there exists a constant $\gamma=\gamma(p,q)>0$
with
\begin{equation}
\phi_{p,q}(0\communique\partial B(0,n))\le e^{-\gamma n}\text{ for large }n.
\end{equation}
\end{prop}

\noindent\textbf{Stochastic comparison}\\ 
Let us first recall the concept of domination for finite measures on a partially
ordered set $E$ . We say that a probability measure $\mu$ dominates a probability measure $\nu$,
if $${\int f\ d\nu}\le{\int f\ d\mu}$$
holds as soon as $f$ in an increasing function.
We also write $\nu\prec\mu$.

If $q'\ge q, q'\ge 1$ and $\frac{p'}{q'(1-p')}\ge\frac{p}{q(1-p)}$, then
$\phi^0_{p,q}\prec\phi^0_{p',q'}$.

\noindent\textbf{Isolation cages}\\ 
Let $F$ be a connected subset of $\Ed$. We say that $F$ is an isolation cage if
$\Ed\backslash F$ as two connected components in $\Ld$. In this case, we denote
by $\Int(F)$ (\resp $\Ext(F)$) the bounded (\resp unbounded) connected component of $\Ed\backslash F$ and  $\Intbarre(F)=\Int(F)\cup F$ (\resp $\Extbarre(F)=\Ext(F)\cup F$).

Let $E$ be a finite connected subset of $\Ed$ and $F$ be an isolation cage. 

Note for $r\in \{0,1\}$, $C_r=\{\forall e\in F;\omega_e=r\}$.

Then, for each $\sigma(\Intbarre(F)\cap E)$ measurable event $M_1$ and each 
each $\sigma(\Extbarre(F)\cap E)$  measurable event $M_2$, we have the decoupling property:
\begin{equation}
\label{decoupling}
\forall (b,r)\in\{0,1\}^2\quad \phi^b_{p,q}(M_1\cap M_2| C_r)=\phi^b_{p,q}(M_1| C_r)\phi^b_{p,q}(M_1\cap M_2| C_r)
\end{equation}
Moreover $$\phi^b_{p,q}(M_1| C_r)=\phi^r_{\Int(F),p,q}(M_1).$$

\subsection{The Potts model}

Let us recall the definition of Gibbs measure in the context of the Potts model. Let $q\ge 2$. We note by $E_q$ a set of cardinal $q$.


For a finite subset $\Lambda$ of $\Zd$, the Hamiltonian on the volume $\Lambda$ is defined by
$$H_{\Lambda}=2\beta\sum_{\substack{e=(x,y)\in\Ed\\ e\cap\Lambda\ne\varnothing}}\1_{\{\omega(x)\ne\omega(y)\}}.$$

Then, we can define for each bounded measurable function $f$ and for each $\omega\in E_q^{\Zd}$,
$${\Pi}_{\Lambda}f({\omega})=\frac1{Z_{\Lambda}({\omega})}
{\miniop{}{\sum}{\eta\in E_q^{\Lambda}}\exp(-H_{\Lambda}(\eta_{\Lambda}{\omega}_{\Lambda^c}))
f(\eta_{\Lambda}{\omega}_{\Lambda^c})}
,$$
where
$$Z_{\Lambda}(\omega)=\sum_{\eta\in E_q^{\Lambda}}\exp(-H_{\Lambda}(\eta_{\Lambda}{\omega}_{\Lambda^c})).$$

For each $\omega$, we will denote by ${\Pi}_{\Lambda}(\omega)$ the 
measure on $E_q^{\Zd}$ which is associated to map $f\mapsto \Pi_{\Lambda} f(\omega)$.
A measure $\mu$ on $E_q^{\Zd}$ is said to be a Gibbs measure for the $q$-state Potts model at inverse temperature $\beta$
when for each bounded measurable function $f$ and each finite subset ${\Lambda}$ of $\Zd$, we have
$$E_{\mu}(f\vert (X_i)_{i\in\Lambda^c})=\Pi_{\Lambda}f\quad\mu\text{ a.s}.$$

For each $r\in E_q$, let us denote by  ${\Pi}_{\Lambda}(r)$ the measure
 ${\Pi}_{\Lambda}(\omega)$ where $\omega$ is the element of $E_q^{\Zd}$ with
$\omega_x=r$ for each $x\in\Zd$. It is known that for each $\beta>0$ and each $q\in\Zd$, the sequence $({\Pi}_{\Lambda}(r))_{\Lambda}$ converges when $\Lambda$ tends to $\Zd$. Let us denote by $\WPt_{q,\beta,r}$ this limit.
By the general theory of Gibbs measures, this limit is necessarily a Gibbs measure -- see for example the reference book by Georgii~\cite{MR89k:82010}.
Häggström, Jonasson and Lyons~\cite{MR2003f:60173} gave a nice characterization of it:
 
\begin{prop} \label{prop24}
Let $q\in \{2,3,\ldots\}$ and $p\in [0,1]$.
Pick a random edge configuration $X\in \{0,1\}^{\Ed}$ according to the random-cluster
measure $\phi^1_{p,q}$. Then, for each {\bf finite} 
connected component $\mathcal{C}$ of $X$ independently, 
pick a spin uniformly from $E_q$, and assign this spin to 
all vertices of $\mathcal{C}$. Finally assign value $r$ to all vertices of
infinite connected components. The $E_q^{\Zd}$-valued random spin
configuration arising from this procedure is then distributed according
to the Gibbs measure $\WPt_{q,\beta, r}$ for the $q$-state Potts model  at inverse temperature $\beta:= -\frac{1}{2}\log(1-p)$.  
\end{prop}

\section{Central Limit Theorem for the random cluster measure}

We begin with a general theorem which gives sufficient conditions for having a central limit theorem for the fluctuations of the size of the intersection of large boxes with the infinite clusters. This will tell us what sort of estimates about random cluster measures can help us. 

\begin{theorem}
\label{general_CLT}
Let $\phi$ be a translation-invariant ergodic measure on $\Sld$.
We suppose that  $\phi$ satisfy to the FKG inequalities and that we have $\theta_{\phi}=\phi(0\communique\infty)>0$.

For each $n\in\Zd$ and $r>0$, let us note the event
$D_{n,r}=\{|C(n)|>r\}$.
Suppose also that there exists a sequence $(r_n)_{n\in\Zd}$ such that the following assumptions together hold:
\begin{itemize}
\item $$(m)\quad \sum_{n\in\Zd} P(+\infty>|C(0)|\ge r_n)<+\infty$$
\item $$(c)\quad \sum_{n\in\Zd} \Cov(D_{0,r_n},D_{n,r_n})<+\infty.$$
\end{itemize}
Then, we have
\begin{itemize}
\item $(S^{*})$ $$\sigma_\phi^2=\sum_{k\in\Zd}\big(\phi_{p,q}(0\communique\infty)\text{ and }k\communique\infty)-\theta_\phi^2\big)<+\infty$$
\item $(CLT)$ $$\frac{|\Lambda_n\cap I|-\theta_{\phi}|\Lambda_n|}{|\Lambda_n|^{1/2}}\Longrightarrow\mathcal{N}(0,\sigma^2_{\phi}).$$
\end{itemize}

\end{theorem}

\begin{proof}
$${\Card{\Lambda_n\cap
I(\omega)}-\theta_{\phi}\Card{\Lambda_n}}=\sum_{k\in\Lambda_n}
f(T^k\omega),$$ where $T^k$ is the translation operator defined by
$T^k(\omega)=(\omega_{k+e})_{e\in\Ed}$ and
$f=\1_{\{\Card{C(0)}=+\infty\}}-\theta_{\phi}$. Moreover, $f$ is an
increasing function and  $\phi$ satisfies the F.K.G. inequalities.
Then, $(f(T^k \omega))_{k\in\Zd})$ is a stationary random field of
square integrable variables satisfying to the F.K.G. inequalities. Therefore,
according to Newman \cite{MR81i:82070}, the Central Limit Theorem
is true if we prove that the quantity
\begin{equation}
\label{laserie} \sum_{k\in\Zd} \Cov(f,f\circ T^k)
\end{equation}
 is finite, which is just proving $(S^{*})$.

Let us define $B=\{|C(0)|=+\infty\}$, and  for each $n\in\Zd$, $A_n=\{|C(n)|=+\infty\}$, 
$\tilde{A}_n=\{|C(n)|\ge r_n\}$ and $\tilde{B}_n=\{|C(0)|\ge r_n\}$
Since $B\subset \tilde{B}_n$ and $A_n\subset \tilde{A}_n$, one has
$\1_{\tilde{B}_n}=\1_{B}+\1_{\tilde{B}_n\backslash B}$ and
$\1_{\tilde{A}_n}=\1_{{A}_n}+\1_{\tilde{A}_n\backslash A_n}$. It follows that
$$\Cov(\1_{\tilde{A}_n},\1_{\tilde{B}_n})-\Cov(\1_{{A}_n},\1_{B})=\Cov(\1_{\tilde{A}_n\backslash A_n},\1_{\tilde{B}_n})+\Cov(\1_{\tilde{B}_n\backslash B},\1_{\tilde{A}_n}),$$
and hence that
\begin{eqnarray*}
|\Cov(\1_{\tilde{A}_n},\1_{\tilde{B}_n})-\Cov(\1_{{A}_n},\1_{B})| & \le & P(\tilde{A}_n\backslash A_n)+P(\tilde{B}_n\backslash B)\\ 
& \le & 2 P(+\infty>|C(0)|\ge r_n)
\end{eqnarray*}

It follows that 
$$\sigma_{\phi}^2\le 2\sum_{n\in\Zd} P(+\infty>|C(0)|\ge r_n)+ \sum_{n\in\Zd} |\Cov(D_{0,r_n},D_{n,r_n})|<+\infty.$$
\end{proof}

The idea of using Newman's theorem to  prove  Central Limit Theorems for the density of infinite clusters in percolation models satisfying to the F.K.G.
inequalities is not new, because it has already been pointed out
by Newman and Schulman \cite{MR83e:82038,MR82i:82047} that
$(S^{*})+FKG\Longrightarrow (CLT)$.
The interest of this theorem is that it splits a problem about infinite clusters into two problems relative to finite clusters:
\begin{itemize}
\item The existence of sufficiently high moments
\item A control of the correlation for the appearance of reasonably large clusters in two points which are separated by a large distance -- note that we can rewrite $\Cov(D_{0,r_n},D_{n,r_n})$ as $\Cov(D_{0,r_n}^c,D_{n,r_n}^c)$.
\end{itemize}

\noindent{\textbf{Example}:} In the case of Bernoulli percolation, the central limit theorem holds for $p>p_c$. 
\begin{proof}
Simply take $r_n=\Vert n\Vert/3$.
The convergence $(m)$ follows for example by the result of Kesten and  Zhang \cite{MR91i:60278}: there
exists $\eta(p)>0$ such that $$\forall n\in\Z_{+}\quad
P(\Card{C(0)}=n)\le\exp(-\eta(p)n^{(d-1)/d}).$$Of course, such a sharp estimate is not
necessary for our purpose. Estimates derived from Chayes, Chayes
and Newman \cite{MR88i:60161}, and from Chayes, Chayes, Grimmett,
Kesten and Schonmann \cite{MR91i:60274} would have been
sufficient. The convergence of $(c)$ is an evidence since $D_{0,r_n}$ and $D_{n,r_n}$ are independent for $\Vert n\Vert>12$.
\end{proof}
\subsection{The case of dimension two}

Let $\Zp=\Z^2+(1/2,1/2)$. For each bond $e=\{a,b\}$ of $\Ldeux$ (\resp $\Ldeuxstar$), let us denote by
$s(e)$ the only subset $\{i,j\}$ of $\Zp$ (\resp $\Z^2$) such that the quadrangle $aibj$ is 
a square. $s$ is clearly an involution. Let us also define $\Ldeuxstar=(\Zp,\Edeuxstar)$, where
$\Edeuxstar=\{s(e);e\in\Edeux\}$.
It is easy to see that  $\Ldeuxstar$ is isomorphic to $\Ldeux$.

For finite $A\subset\Zp$, we denote by $\peierls(A)$  the Peierls contour associated 
to $A$, that is

$$\peierls(A)=\{e\in\Edeux;\1_A\text{is not constant on }s(e)\}.$$
There exists a finite family of cycles and paths -- the so-called Peierls contours -- such that $\peierls(A)$ is the set of vertices visited by  them.
It is known that, provided that $A$ is a bounded connected subset of $\Ldeuxstar$, there exists a unique cycle $\Gamma(A)$ which is a Peierls contour and surrounds $A$.

Note that in the two-dimensional lattice, the set of bonds of a cycle
forms an isolation cage.
So if $\gamma$ is a cycle, we will simply denote by $\Intbarre(\gamma)$ the
set   $\Intbarre(E)$, where $E$ is the set of bonds of $\gamma$. 

Consider now the map
\begin{eqnarray*}
\{0,1\}^{\Edeux} & \to & \{0,1\}^{\Edeuxstar}\\
\omega & \mapsto & \omega^{*}=(1-\omega_{s(e)})_{e\in\Edeuxstar}
\end{eqnarray*}
For $\eta\in\{0,1\}^{\Edeuxstar}$, we also denote by $\eta^{*}$ the only $\omega\in\{0,1\}^{\Edeux}$ such that $\omega^{*}=\eta$.

For each subset $A$ of $\{0,1\}^{\Edeux}$ (\resp  $\{0,1\}^{\Edeuxstar}$), we denote by $A^{*}$ the set $\{\omega^{*};\omega\in A\}$.

The following planar duality between planar random cluster measures is now well known: 
let us define $p^{*}$ to be the 
unique element of $[0,1]$ which satisfies to $F(p)F(p^{*})=1$, with
$F(x)=\frac1{\sqrt{q}}\frac{x}{1-x}$.
and also define a map $t$ by
\begin{eqnarray*}
\{0,1\}^{\Edeux} &\to & \{0,1\}^{\Edeuxstar}\\
\omega & \mapsto & (\omega_{e+(1/2,1/2)})_{e\in\Edeuxstar}
\end{eqnarray*}
Then, for each $p\in [0,1]$, $b\in\{0,1\}$ and each event $A$, we have
$$\phi^{b}_{p,q}(A)=t\phi^{1-b}_{p^{*},q}(A^{*}).$$

Let us note $r(q)=p_g(q)^{*}$. Since $p_g(q)>0$, we have $r(q)<1$.
Note that it is believed that $p_g(q)=p_c(q)$. 
As was noted by Grimmett and Piza~\cite{MR99a:60115}, the fact that  $p_g(q)=p_c(q)$ would imply that $p_c$ is the  solution of the equation $x=x^{*}$, \ie $p_c=\frac{\sqrt q}{1+\sqrt q}$.
So, it follows that we have $r(q)=p_c(q)$ provided that $p_g(q)=p_c(q)$.
When $d=2$, this widely believed conjecture has already be proved
for $q=1,2$ and $q\ge 22$ -- see the Saint-Flour notes by Grimmett~\cite{MR99c:60216}.   

\begin{lemme}
\label{underive}
Let $d=2$ and $p<p_g(q)$. There exists $K\in (0,+\infty)$ with
$$\forall n\in\N\quad \phi_{p,q}(|C(0)|\ge n)\le K\exp(-\gamma(p,q) \sqrt{n/2}).$$
\end{lemme}
\begin{proof}
Suppose that $n\ge 16$ and denote by $r$ the integer part of $\sqrt{n/2}-1$.
Let $T=\{k\in\N; C(0)\cap \partial B(0,k)\ne\varnothing\}$ and $R=\max T$.\\
Suppose $|C(0)|\ge n$:  we have $C(0)\subset B(0,R)$, so $$n\le |C(0)|\le |B(0,R)|=1+2R(R+1).$$ It follows that $r\le R$. Since $C(0)$ is connected, we have $0\communique \partial B(0,r)$.
Then $\phi_{p,q}(|C(0)|\ge n)\le \phi_{p,q}(0\communique \partial B(0,r))$.
The result follows then from Proposition~\ref{gr-pi}.
\end{proof}

When $d=2$, it is known that $p_g\le p_c\le\frac
{\sqrt{q}}{1+\sqrt q}$. It follows that 
$p_g^{*}\ge p_c^{*}\ge\frac{\sqrt{q}}{1+\sqrt q}$. Then $(p>p_g^{*}) \Longrightarrow (p\ne p_c)$. Since it is known that $\mathcal{R}_{p,q}$ is a single when $p\ne p_c$ and $d=2$, it follows that there is an unique
random cluster measure for $p>r(q)$. Then, we simply write $\phi_{p,q}$ without any superscript.

\begin{lemme}
\label{dimdeux}
We suppose here that $d=2$.  For each $p>r(q)$,
there exists $K\in (0,+\infty)$ with 
$$\forall n\ge 1\quad\phi_{p,q}(|C(0)|=n)\le K n e^{-\gamma(p^{*},q) \sqrt{n}}.$$
Note that $\gamma(p^{*},q)>0$.
\end{lemme}
\begin{proof}
We use here a duality argument. 
Let $p>r(q)$ and note $A=\{|C(0)|=n\}$. We have $\phi_{p,q}(A)=t\phi_{p^{*},q}(A^{*})$. In this case 
$$t^{-1}(A^{*})=\big\{\text{
\begin{tabular}{c}
there exists at least one open cycle surrounding $(0,1/2)$\\
Those of these cycles which minimizes the distance to $(0,1/2)$\\ 
surrounds exactly $n$ closed bonds. 
\end{tabular}
}\big\}.$$
The number of bounds used by this cycle is at least $2n+2$.
Moreover, the position of the first intersection of this cycle with the positive $x$-axis is at most $n$.
So 
$$t^{-1}(A^{*})\subset\miniop{n}{\cup}{k=1}\{|C(ke_1)|\ge 2n\}.$$
It follows then from lemma~\ref{underive} that

$$\phi_p(A)=\phi_{p^{*},q}(A^{*})\le K ne^{-\gamma(p^{*},q)\sqrt{n}}.$$

\end{proof}

The goal of the next lemma is to bound the covariance of two decreasing events who are defined by the state of the bonds in two boxes separated by a large distance. It is clear that it does not pretend to originality and that its use could have been replaced by those of an analogous result of the literature, \eg 
 Theorem 3.4 of Alexander~\cite{MR99e:60211} joined to its Remark 3.5.
Nevertheless, we preferred to present our lemma because its proof is rather short and allows an instructive comparison with the case of  an higher dimension
which will be studied after. 

\begin{lemme}
\label{covdeux}
Let $q\ge 1$. For each $p>r(q)$, there exists $C>0$ and $\alpha>0$ such that for each couple of boxes $\Lambda_1$ and $\Lambda_2$ and each pair of monotone
events $A$ and $B$, with $A$ (\resp $B$) $\sigma(\Lambda_1)$ (\resp  $\sigma(\Lambda_2)$ ) measurable, we have
$$|\Cov_{\phi}(A,B)|=| \phi_{p,q}(A\cap B)-\phi_{p,q}(A)\phi_{p,q}(B)|\le C|\partial\Lambda|e^{-\alpha d(\Lambda_1,\Lambda_2)}.$$
\end{lemme}
\begin{proof}
Since $\Cov_{\phi}(A,B)=-\Cov_{\phi}(A^c,B)=\Cov_{\phi}(A^c,B^c)=-\Cov_{\phi(A,B^c)}$, we can assume that $A$ and $B$ are decreasing events.
We can also assume without loss of generality that $\Lambda_1=\{-n,\dots,n\}\times\{-p,\dots,p\}$. Put $\Lambda_1^{*}=\{-n+1/2,\dots, n-1/2\}\times\{-p+1/2,\dots, p-1/2\}$.
For $x\in\Zp$ and $\omega\in\Omega$, let us define $C^{*}(x)$ to be the connected component of $x$ in the configuration $\omega^{*}$.  
Let now be $F$  the random set defined by
$$F=\Gamma(\miniop{}{\cup}{y\in\Lambda_1^{*}}C^{*}(y))$$
and note $V=\{\Intbarre(F)\cap\Lambda_2=\varnothing\}$.
The following facts are elementary, but relevant:
\begin{itemize}
\item For every curve $\gamma$ surrounding the origin,  the event 
$\{F=\gamma\}$ is $\sigma({\Intbarre(\gamma)})$-measurable. 
\item For any subset $\gamma$ of $\Zd$, $\{F=\gamma\}\subset W_{\gamma}=\{\forall e\in \gamma;\omega_e=1\}$.
\end{itemize} 
Remember that if $T$ is a $\sigma({\Extbarre(\gamma)})$-measurable event and
$R$ a $\sigma(\Intbarre(\gamma))$-measurable event, we have
the following decoupling property:
\begin{eqnarray*}
\phi_{p,q}(R\cap T|W_{\gamma}) & = & \phi_{p,q}(R|W_{\gamma})\phi_{p,q}(T|W_{\gamma})\\
\end{eqnarray*}
So, if $\Intbarre(\gamma)\cap\Lambda_2=\varnothing$ we have

\begin{eqnarray*}
\phi_{p,q}(A\cap B\cap \{F=\gamma\}) & = & \phi_{p,q}(A\cap B\cap \{F=\gamma\}\cap W_{\gamma})\\
 & = & \phi_{p,q}(W_\gamma)\phi_{p,q}(A\cap \{F=\gamma\})\cap B | W_{\gamma})\\
& = &  \phi_{p,q}(W_\gamma) \phi_{p,q}(A\cap \{F=\gamma\})| W_{\gamma})
\phi_{p,q}(B | W_{\gamma})\\
& = & \phi_{p,q}(A\cap \{F=\gamma\})\cap W_{\gamma})
\phi_{p,q}(B | W_{\gamma})\\
& = & \phi_{p,q}(A\cap \{F=\gamma\})
\phi_{p,q}(B | W_{\gamma})\\
& \le & \phi_{p,q}(A\cap \{F=\gamma\})\phi_{p,q}(B)
\end{eqnarray*}

If we sum over suitable values of $\gamma$, we get
$$\phi_{p,q}(A\cap B\cap V)\le  \phi_{p,q}(A\cap V)\phi_{p,q}(B)\le  \phi_{p,q}(A)\phi_{p,q}(B).$$
Since $A$ and $B$ are decreasing events, they are positively correlated, then
$$0\le\phi_{p,q}(A\cap B)-\phi_{p,q}(A)\phi_{p,q}(B)\le\phi_{p,q}(A\cap B\cap V^c)\le\phi_{p,q}(V^c).$$
Since $$V^c\subset\miniop{}{\cup}{y\in\partial\Lambda^{*}}\{y\communique\partial B(y,d(\Lambda_1,\Lambda_2))\},$$
the result follows from the inequality of Grimmett and Piza.
\end{proof}

\subsection{The case of general $d$}
The goal of the next lemma is also to bound the covariance of two positive events who are defined by the state of the bonds in two boxes separated by a large distance.

Unlike  the proof of lemma~\ref{covdeux}, it can not  use duality 
arguments. We nevertheless attempt to present this proof in a form which is as close as possible of those of lemma~\ref{covdeux} to highlight the differences and the similarities between them.


Not that it is somewhat inspired by the proof of Grimmett~\cite{MR97b:60169} for the uniqueness of the random-cluster when $p$ is large. 

\begin{lemme}
\label{covar}
Let $q\ge 1$. For each $p>\frac{q}{q+4^{-d}}$, there
exists $\alpha(p,q)>0$ such that
for each couple of finite connected volumes $\Lambda_1$ and $\Lambda_2$ and each pair of monotone
events $A$ and $B$, with $A$ (\resp $B$) $\sigma(\Lambda_1)$ (\resp  $\sigma(\Lambda_2)$ ) measurable, we have
$$|\Cov_{\phi}(A,B)|=| \phi_{p,q}(A\cap B)-\phi_{p,q}(A)\phi_{p,q}(B)|\le C|\partial\Lambda|e^{-\alpha d(\Lambda_1,\Lambda_2)}.$$
\end{lemme}

\begin{proof}

We begin by a topological remark: let $A$ be  a finite connected subset of 
$\Zd$ and consider the $*$-connected components of $\partial A$: it is not difficult to see that there is exactly one of the connected components, say $B$, which surrounds $A$. Let us also define $W(A)=\{e\in\Ed; e\cap B\ne\varnothing\}$. Since $B$ is $*$-connected, $W(A)$ is connected in $\Ld$; $W(A)$ also surrounds $A$. Note that $W(A)$ is the analogous to a surrounding Peierls contour in the two dimensional lattice.

Now suppose as previously that $A$ and $B$ are decreasing events.
Given a configuration $\omega$, say that a point $x\in\Zd$ 
is wired if each of the $2^d$ bonds attached to $x$ satisfy to $\omega_e=1$.
Otherwise, say that $x$ is free.
Let us define $D(\omega)$ to be the set of points in $\Zd\backslash\Lambda_1$ which can be connected to $\Lambda_1$ using only free vertices -- the origin of the path in $\Lambda_1$ does not need to be free.
By definition of $D$, $(\Lambda_1 \cup D)$ is a connected set. 

Let us note $F=W(D)$ and $V=\{ \Lambda_2\cap D=\varnothing\}$.

Note that $V$ is an increasing event.

The following facts are elementary, but relevant:
\begin{itemize}
\item For every isolation cage $\gamma$ surrounding $\Lambda_1$,  the event 
$\{F=\gamma\}$ is $\sigma({\Intbarre(\gamma)})$-measurable. 
\item For any subset $\gamma$ of $\Ed$, $\{F=\gamma\}\subset W_{\gamma}=\{\forall e\in \gamma;\omega_e=1\}$.
\end{itemize} 
Remember that if $T$ is a $\sigma({\Extbarre(\gamma)})$-measurable event and
$R$ a $\sigma(\Intbarre(\gamma))$-measurable event, we have
the following decoupling property:
\begin{eqnarray*}
\phi^0_{p,q}(R\cap T|W_{\gamma}) & = & \phi^0_{p,q}(R|W_{\gamma})\phi^0_{p,q}(T|W_{\gamma})\\
\end{eqnarray*}
So, if $\gamma$ does not touch nor surround $\Lambda_2$,
we have

\begin{eqnarray*}
\phi^0_{p,q}(A\cap B\cap \{F=\gamma\}) & = & \phi^0_{p,q}(A\cap B\cap \{F=\gamma\}\cap W_{\gamma})\\
 & = & \phi^0_{p,q}(W_\gamma)\phi^0_{p,q}(A\cap \{F=\gamma\})\cap B | W_{\gamma})\\
& = &  \phi^0_{p,q}(W_\gamma) \phi^0_{p,q}(A\cap \{F=\gamma\})| W_{\gamma})
\phi^0_{p,q}(B | W_{\gamma})\\
& = & \phi^0_{p,q}(A\cap \{F=\gamma\})\cap W_{\gamma})
\phi^0_{p,q}(B | W_{\gamma})\\
& = & \phi^0_{p,q}(A\cap \{F=\gamma\})
\phi^0_{ p,q}(B | W_{\gamma})\\
& \le & \phi^0_{p,q}(A\cap \{F=\gamma\})\phi^0_{p,q}(B)\\
\end{eqnarray*}

If we sum over suitable values of $\gamma$, we get
$$\phi^0_{p,q}(A\cap B\cap V)\le  \phi^0_{p,q}(A\cap V)\phi^0_{p,q}(B)\le  \phi^0_{p,q}(A)\phi^0_{p,q}(B).$$
Since $A$ and $B$ are decreasing events, they are positively correlated, then
$$0\le\phi^0_{p,q}(A\cap B)-\phi^0_{p,q}(A)\phi^0_{p,q}(B)\le\phi^0_{p,q}(A\cap B\cap V^c)\le\phi^0_{p,q}(V^c).$$

Since $V$ is an increasing event, we can use the stochastic domination
of a product measure by  $\phi^0_{p,q}$:
$\phi^0_{r,1}\prec \phi^0_{p,q}$, with 
$r=p/(p+(1-p)q)$, so 
$\phi^0_{p,q}(F^c)\le \phi^0_{r,1}(F^c)$.
Now, a Peierls-like counting argument gives
$$\phi^0_{r,1}(F^c)
 \le |\partial\Lambda_1| (2^d-1) (2^d(2^d-1))^n (1-r)^n,$$ where $n=d(\Lambda_1,\Lambda_2)-1,$
which completes the proof. 
\end{proof}

\begin{theorem}
\label{CLT} For each $q\ge 1$, there exists $r(q)<1$ such that, for $p>r(q)$, $\mathcal{R}_{p,q}$ consists in an unique measure $\phi_{p,q}$. Moreover, if we note $\theta(p,q)=\phi_{p,q}(0\in I)$, we have
$$\frac{\Card{\Lambda_n\cap
I}-\theta(p,q)\Card{\Lambda_n}}{{\Card{\Lambda_n}^{1/2}}}\Longrightarrow\mathcal{N}(0,\sigma_p^2),$$
where $I$ is the infinite cluster for FK percolation and
$$\sigma_{p,q}^2=\sum_{k\in\Zd}\big(\phi_{p,q}(0\communique\infty\text{ and }k\communique\infty)-\theta(p,q)^2\big).$$
Note also that we can take $r(1)=p_c(\Zd)$, and $r(q)=p_g^{*}$ for $d=2$.
\end{theorem}

\begin{proof}
The uniqueness of the random cluster measure for $p$ close to 1 as be proved
by Grimmett~\cite{MR97b:60169}. So it only remains to prove that it holds
for $p>r(q)$ in the cases where we have announced a convenient value for
$r(q)$. When $q=1$, the uniqueness is obvious and we have already remarked
that there was uniqueness for $d=2$ and $p> p_g^{*}$.

Let us now prove the Central Limit Theorem.
We will apply Theorem~\ref{general_CLT} to the sequence $r_n=\Vert n\Vert/4$.

Let us show that
$$\sum_{n\in\Zd} P(+\infty>|C(0)|\ge \frac{\Vert n\Vert}4)$$ converges.
\begin{itemize}
\item For $d\ge 3$ and $p$ sufficiently close to 1, this follows from the estimate of Pisztora~\cite{MR97d:82016}: for each $b\in\{0,1\}$ and each $p$ which is sufficiently close to 1, the exist a constant $a=a(p,q)$ with
$$
\forall n\ge 0\quad\phi_{p,q}^b(|C|=n) \le \exp\big(-a n^{(d-1)/d}\big).
$$
\item For $d=2$ and $p>p_{g}(q)^{*}$, it follows from our lemma~\ref{dimdeux}.
\end{itemize}
Now, it remains to prove that
\begin{equation}
\sum_{n\in\Zd} \Cov(\1_{\tilde{A}_n},\1_{\tilde{B}_n})<+\infty,
\end{equation}
with $\tilde{A}_n=\{|C(0)|\ge r_n\}$ and $\tilde{B}_n=\{|C(n)|\ge r_n\}$.

Put $\Lambda_n=B(n,\Vert n\Vert/3)$ and $\Lambda'_n=B(0,\Vert n\Vert/3)$.
It is clear that $\tilde{A}_n$ ( \resp $\tilde{B}_n$) is
$\sigma({\Lambda_n})$ (\resp $\sigma({\Lambda'_n})$) measurable.
It is obvious that  $\tilde{A}_n$ and  $\tilde{B}_n$ are increasing events.
Then, we can apply lemma \ref{covar}. Since $d(\Lambda_n,\Lambda'_n)\ge n/3$, it follows that
$$0\le \Cov(\1_{\tilde{A}_n},\1_{\tilde{B}_n})\le K n^{d-1}e^{-\frac{\alpha(p,q)}3n},$$
which forms a convergent series as soon as $p>p(q)$. 
When $d=2$, the result follows similarly from lemma~\ref{covdeux}.
\end{proof}

\section{Coloring of random clusters}

If $G$ is a subgraph of $\Ld$, $r\in\R$ and if $\nu$ is a probability
measure on $\R$, we will define the color-probability $P^{G,\nu,r}$
as follows: $P^{G,\nu,r}$ is the unique measure on
$(\R^{\Zd},\mathcal{B}(\R^{\Zd}))$ under which the canonical
projections $X_i$ -- defined, as usual by $X_i(\omega)=\omega_i$
-- satisfy
\begin{itemize}
\item For each $i\in\Zd$, the law of $X_i$ is 
\begin{itemize} \item $\nu$ if $|C(i)|<+\infty$.
\item $\delta_{r}$ otherwise.
\end{itemize}
\item For each independent set $S\subset\Zd$, the variables
$(X_i)_{i\in S}$ are independent.
\item For each connected set $S\subset\Zd$, the variables
$(X_i)_{i\in S}$ are identical.
\end{itemize}

Let $\phi$ be a measure on $\Sld$.

The randomized color-measure associated to $\phi$ is defined by 
$$P^{\phi,\nu,r}=\int_{\Sld}
P^{G,\nu,r}\ d\phi(G).$$

Of course, our aim is to specialize $\phi$ to a Fortuin-Kastelein measure
$\phi=\phi^b_{p,q}$, with $b\in\{0,1\}$, but the reasoning that we will made only depend from the existence of
\begin{itemize}
\item stationarity and ergodicity of $\phi$
\item moments condition for the size of finite clusters
\item existence of a central limit theorem for size of the intersection of infinite clusters with large boxes
\end{itemize} 
In order to make easier the later use of these results according to the progress
that would be made in random-cluster or related models, we try to expose
our results in the most possible generality, and then to apply them according to the results that we have nowadays.

To motivate this work, let us give some examples of models covered by randomized color-measures when $\phi=\phi^b_{p,q}$.
\begin{itemize}
\item 
The case $q=1$ is a generalization of the divide and color model of Häggström~\cite{MR2003i:60172}, which has already been studied in a earlier paper of the author~\cite{MR2003g:60166}. 
\item The most celebrated of the randomized color-measure is obtained 
when $q\ge 2$ is an integer and $\nu=\frac1q(\delta_1+\delta_2+\dots \delta_q)$. In this case $P^{\phi,\nu,r}$ is the Gibbs measure
$\WPt_{q,\beta, r}$ for the $q$-state Potts model on
$\Zd$ at inverse temperature $\beta:= -\frac{1}{2}\log(1-p)$, according
to Proposition~\ref{prop24}.
It includes of course the case of the Ising model.
\item If $n_1,n_2,\dots n_s$ are positive integers with $n_1+n_2+\dots n_s=q$
and we take $\nu=\frac1q(n_1\delta{1}+n_2\delta_{2}+\dots n_s\delta_s)$, we
obtain a fuzzy Potts model. It obviously follows from the previous example and
the definition of a fuzzy Potts model.
\end{itemize}

In all this section, we will suppose that $\phi$ satisfy to the following assumptions:
\begin{itemize}
\item (E): $\phi$ is a translation-invariant ergodic measure on $\Sld$.
\item (M): $\exists \alpha>d\quad \sum_{k=1}^{+\infty} k^{\alpha} P(\Card{C(0)}=k)<+\infty.$
\end{itemize}

When $\theta_{\phi}=\phi(0\communique\infty)>0$, the following assumption will also be considered:
$$(CLT): \exists \sigma^2_{\phi}>0, \frac{|\Lambda_n\cap I|-\theta_{\phi}|\Lambda_n|}{|\Lambda_n|^{1/2}}\Longrightarrow\mathcal{N}(0,\sigma^2_{\phi}).$$

We begin with a general property of randomized color-measures.
\begin{theorem}
\label{ergodicite}
$P^{\phi,\nu,z}$ is translation invariant and the action of $\Zd$ on $P^{\phi,\nu,z}$ is ergodic.
\end{theorem}
\begin{proof}
Let $\Omega=\{0,1\}^{\Ed}$,   $\Omega_t=[0,1]^{\Zd}$ , $\Omega_c=\R^{\Zd}$,
$\Omega_3=\Omega\times \Omega_t\times \Omega_c$ and consider the action of $\Zd$ on
$(\Omega_3,\bor[\Omega_3],\phi\otimes U[0,1]^{\otimes\Zd}\otimes \nu^{\otimes\Zd})$. Since $\phi\otimes U[0,1]^{\otimes\Zd}\otimes \nu^{\otimes\Zd}$ is a direct product of an ergodic measure by two mixing measures, it follows that the action of $\Zd$ on $\Omega_3$ is ergodic -- see Brown~\cite{Brown} for instance.

Let us define $f:\Omega_3\to\R$ by
$f(\omega,\omega_t,\omega_c)=z$ if 
$|C(0)(\omega)|=+\infty$ or 
if exists $x,y\in C(0)(\omega)$, with $x\ne y$ and $\omega_t(y)=\omega_t(x)$.
Otherwise, we define $f(\omega,\omega_t,\omega_c)$ to be equal 
to $\omega_c(x)$, where $x$ is the unique element of $C(0)(\omega)$ such that
$(\omega_t)(x)=\max\{\omega_t(y);y\in C(0)\}$.
Now, if we define $(X_k)_{k\in\Zd}$, by $X_k(\omega,\omega_t,\omega_c)=
f(\theta_k(\omega,\omega_t,\omega_c))$, it is not difficult that the law
of $(X_k)_{k\in\Zd}$ under $\phi\otimes U[0,1]^{\otimes\Zd}\otimes \nu^{\otimes\Zd}$ is $P^{\phi,\nu,r}$.
Since a factor of an ergodic system is an ergodic system -- see also Brown~\cite{Brown} -- , it follows that $P^{\phi,\nu,z}$ is ergodic under the action of $\Zd$.
\end{proof}

\subsection{Normal fluctuations of sums for color-measures}

\begin{theorem}
\label{quenched}
Suppose that $\phi$ satisfy to (E) and (M) and
let $\nu$ be a probability measure on $\R$ with a
second moment. We put $m=\int_{\R} x\ d\nu(x)$ and
$\sigma^2=\int_{\R} (x-m)^2\ d\nu(x)$.

 For $\phi$-almost all $G$, the following holds:

 $$\frac1{\Card{\Lambda_n}^{1/2}}\big(\sum_{x\in
\Lambda_n\backslash I }(X(x)-m)\big)\Longrightarrow\mathcal{N}(0,
\chi^f(\phi)\sigma^2)$$ where 
$I(G)=\{x\in\Zd; x\communique\infty\}$.
\end{theorem}

The following lemma will be very useful.

\begin{lemme}
\label{quoca} For each subgraph $G$ of $\Ld$, let us denote by
$(A_i)_{i\in J}$ the partition of $G$ into connected components.

Suppose that $\phi$ satisfy to (E) and (M).
 Then, we have
for $\phi$-almost all $G$:
$$\lim_{n\to\infty}\frac1{\Card {\Lambda_n}}\sum_{i\in
J;\Card{A_i}<+\infty}\Card{A_i\cap \Lambda_n}^2=\chi^f(\phi),$$

where

$$\chi^f(\phi)=\sum_{k=1}^{+\infty} k \phi(\Card{C(0)}=k).$$
\end{lemme}
\begin{proof}
Let us define $C'(x)$ by

\begin{equation*}
C'(x)= \begin{cases} C(x) & \text{ if }\Card{C(x)}<+\infty\\
\varnothing & \text{otherwise}
\end{cases}
\end{equation*}
 and $C'_n(x)=C'(x)\cap \Lambda_n$.

It is easy to see that
\begin{equation}
\label{remarquable}
 \sum_{i\in J;\Card{A_i}<+\infty}\Card{A_i\cap
\Lambda_n}^2=\sum_{x\in \Lambda_n}\Card{C'_n(x)}.
\end{equation}
 We have
$\Card{C'_n(x)}\le\Card{C(x)}$, and the equality holds if and only if
$C'(x)\subset \Lambda_n$.

We have
$$\sum_{k=1}^{+\infty} n^{d-1}\phi(+\infty|C(0)|\ge n^{d/\alpha})\le\int_0^{+\infty} x^{d-1}\phi(+\infty>|C(0)|\ge n^{d/\alpha})=\frac1d
\int_{{|C(0)|<+\infty}}
 |C(0)|^{\alpha}\ d\phi.$$
 It follows from a standard Borel-Cantelli argument
that for $\mu_p$-almost all $G$, there exists a (random) $N$ such
that $$\forall n\ge N\quad \max_{x\in \Lambda_n} \Card{C'(x)}\le
n^{d/\alpha}.$$ If follows that for each $x\in
\Lambda_{n-n^{d/\alpha}}$, $C'(x)$ is completely inside
$\Lambda_n$, and therefore $C'(x)=C'_n(x).$ Then,
\begin{eqnarray*}
\sum_{x\in \Lambda_{n-n^{d/\alpha}}}\Card{C'(x)} & \le &
\sum_{x\in \Lambda_n}\Card{C'_n(x)}\le \sum_{x\in
\Lambda_n}\Card{C'(x)}\\ \frac1{\Card{\Lambda_n}}\sum_{x\in
\Lambda_{n-n^{d/\alpha}}}\Card{C'(x)} & \le &
\frac1{\Card{\Lambda_n}}\sum_{x\in \Lambda_n}\Card{C'_n(x)}\le
\frac1{\Card{\Lambda_n}}\sum_{x\in \Lambda_n}\Card{C'(x)}\\
\end{eqnarray*}

By the ergodic Theorem, we have $\phi$-almost surely: $$\lim_{n\to
+\infty} \frac1{\Card{\Lambda_n}}\sum_{x\in
\Lambda_n}\Card{C'(x)}=\int \Card{C'(0)}\ d\phi=\chi^f(\phi).$$ Since
$\lim_{n\to +\infty}\frac{\Card{\Lambda_{n-n^{d/\alpha}}}}{\Card{\Lambda_n}}=1$ , the result follows.

\end{proof}

\begin{proof}
$$\sum_{x\in \Lambda_n\backslash I
}(X(x)-m)=\miniop{+\infty}{\sum}{i=1}\Card{C'_n(a_i)}(X(a_i)-m)$$
Then $$\frac1{\Card{\Lambda_n}^{1/2}}\sum_{x\in
\Lambda_n\backslash I
}(X(x)-m)=\big(\frac{s_n^2}{\Card{\Lambda_n}}\big)^{1/2}\frac1{s_n}\miniop{+\infty}{\sum}{i=1}\Card{C'_n(a_i)}(X(a_i)-m),$$
with $$s_n^2=\miniop{+\infty}{\sum}{i=1}\Card{C'_n(a_i)}^2.$$

By lemma \ref{quoca}, we have for $\phi$-almost all $G$
$\lim_{n\to +\infty}\frac{s_n^2}{\Card{\Lambda_n}}=\chi^f(\phi)$.

Now, it remains to prove that
\begin{equation}
\frac1{s_n}\miniop{+\infty}{\sum}{i=1}\Card{C'_n(a_i)}(X(a_i)-m)\Longrightarrow\mathcal{N}(0,\sigma^2).
\end{equation}
Therefore, we will prove that for $\phi$-almost all $G$, the
sequence $Y_{n,k}=\Card{C'_n(a_i)}(X(a_i)-m)$ satisfies the
Lindeberg condition. For each $\epsilon>0$, we have
\begin{eqnarray*}
\sum_{k=1}^{+\infty}\frac1{s_n^2}\int_{\vert Y_{n,k}\vert\ge
\epsilon s_n} Y_{n,k}^2\ d P^{G,\nu} & = &
\sum_{k=1}^{+\infty}\frac{\Card{C'_n(a_k)}^2}{s_n^2}\int_{
\Card{C'_n(a_k)}\vert x\vert\ge \epsilon s_n} (x-m)^2\ d\nu(x)\\ &
\le & \int_{\vert x\vert\ge \frac{\epsilon}{\eta_n}} (x-m)^2\
d\nu(x),
\end{eqnarray*}
with $\eta_n=\frac{\sup_{k\ge 1}\Card{C'_n(a_k)}}{s_n}$. Thus, the
Lindeberg condition is fulfilled if $\lim \eta_n=0$. But we have
already seen that $s_n\sim (\chi^f(\phi)\Card{\Lambda_n})^{1/2}$,
whereas $\sup_{k\ge 1}\Card{C'_n(a_k)}=O(n^{d/\alpha})=o(n^{d/2})$. This
concludes the proof.

\end{proof}

\begin{theorem}
\label{annealed} 
Let $\phi$ be a measure on $\Sld$ that satisfy to $(E)$ and $(M)$.
Let $\nu$ be a probability measure on $\R$ with a
second moment. We put $m=\int_{\R} x\ d\nu(x)$ and
$\sigma^2=\int_{\R} (x-m)^2\ d\nu(x)$. 
Let also $z\in\R$.

\begin{itemize}
\item If $\theta_{\phi}=0$, then we have under $P^{\phi,\nu,z}$
 $$\frac1{\Card{\Lambda_n}^{1/2}}\big(\sum_{x\in
\Lambda_n\backslash I }(X(x)-m)\big)\Longrightarrow\mathcal{N}(0,
\chi^f(\phi)\sigma^2)$$
\item If $(\theta_{\phi}>0)$ and $(CLT)$ hold , then we have under $P^{\phi,\nu,z}$
$$\frac1 {\Card{\Lambda_n}^{1/2}}\big(\sum_{x\in
\Lambda_n}X(x)-((1-\theta_{\phi})m+\theta(p)r)\Card{\Lambda_n})\big)\Longrightarrow \mathcal{N}(0,\chi^f(\phi)\sigma^2+(z-m)^2\sigma_{\phi}^2).$$
\end{itemize}
\end{theorem}

\begin{proof}
In this proof, it will be useful to consider $G$ as a random variable.
Let $\Omega=\Sld\times\R^{\Zd}$ and define the probability $\mathbb{P}$ on
$\mathcal{B}(\Omega)$ as a skew-product: for measurable 
$A\times B\in\mathcal{B}(\Sld)\times \mathcal{B}(\R^{\Zd})$, we have 
$\mathbb{P}(A\times B)=\int_A P^{G,\nu,z}(B)\ d\phi(G)$. 
Then, the law of the marginals $G$ and $X$ are $\mathbb{P}_G=\phi$ and 
$\mathbb{P}_X=P^{\phi,\nu,z}$.
Ase usually, the letter $\E$ will be used to denote an expectation -- or a conditional expectation -- under $\P$. 
  
Rearranging the terms of the sum, we easily obtain
$$\big(\sum_{x\in
\Lambda_n}X(x)-((1-\theta_{\phi})m+\theta_{\phi}z)\Card{\Lambda_n})\big)=\sum_{x\in
\Lambda_n\backslash I
}(X(x)-m)+(z-m)(\Card{I\cap\Lambda_n}-\Card{\Lambda_n}\theta_{\phi})$$
We will now put $$Q_n=\frac1
{\Card{\Lambda_n}^{1/2}}\big(\sum_{x\in
\Lambda_n}X(x)-((1-\theta_{\phi})m+\theta_{\phi}z)\Card{\Lambda_n})\big),$$
and define $$\forall t\in\R\quad \phi_{n,z}(t)=\E \exp(iQ_nt).$$

Thereby, we have $$\phi_{n,z}(t)=\E \exp(-
\frac{it}{\Card{\Lambda_n}^{1/2}}\sum_{x\in \Lambda_n\backslash I
}(X(x)-m)+(z-m)(\Card{I\cap\Lambda_n}-\Card{\Lambda_n}\theta_{\phi}))$$
Conditioning by $\sigma(G)$ and using the fact that $I$ is
$\sigma(G)$-measurable, we get $\phi_{n,z}(t)=\E f_n(t,.)
g_{n}((z-m)t,.)$, with
\begin{eqnarray*}
f_n(t,\omega)& = & \E \exp(-
\frac{it}{\Card{\Lambda_n}^{1/2}}\sum_{x\in \Lambda_n\backslash I
}(X(x)-m)\vert\sigma(G)\\ & = & \int \exp(-
\frac{it}{\Card{\Lambda_n}^{1/2}}\sum_{x\in \Lambda_n\backslash
I(\omega) }(X(x)-m)\ dP^{G(\omega),\nu} \end{eqnarray*}
 and
\begin{eqnarray*}
g_{n}(t,\omega)& = &  \exp(-
\frac{it}{\Card{\Lambda_n}^{1/2}}(\Card{I(\omega)\cap\Lambda_n}-\Card{\Lambda_n}\theta_{\phi}))
\\
\end{eqnarray*}
By Theorem~\ref{quenched} we have for each $t\in\R$ and
$P^{\phi,\nu,z}$-almost all $\omega$: $\miniop{}{\lim}{n\to + \infty}
f_{n}(t,\omega) =\exp(-\frac{t^2}{2}\chi^f(\phi)\sigma^2)$ Then, by
dominated convergence $$\lim_{n\to +\infty}\E
(f_n(t,.)-\exp(-\frac{t^2}2\chi^f(\phi)\sigma^2))
g_{n}((z-m)t,.)=0.$$ Next
\begin{eqnarray*}
\lim_{n\to +\infty}\E f_n(t,.) g_{n}((z-m)t,.) & = & \lim_{n\to
+\infty}\exp(-\frac{t^2}2\chi^f(\phi)\sigma^2)\E g_{n}((z-m)t,.)\\ &
= &
\exp(-\frac{t^2}2\chi^f(\phi)\sigma^2)\exp(-\frac{t^2}2(z-m)^2\sigma^2_p)
\end{eqnarray*}
where the last equality follows from Proposition~\ref{CLT}. We
have just proved that $$\miniop{}{\lim}{n\to\infty}
\phi_{n,z}(t)=\exp(-\frac{t^2}2(\chi^f(\phi)\sigma^2+(z-m)^2\sigma^2_p)).$$
The result follows from the Theorem of Levy.
\end{proof}

\subsection{Fluctuation of the empirical vector associated to coloring models}
\hspace{0cm}\\

 \textbf{Definition} Let $q$ be an integer with $q\ge 2$. For every $r\in\{1,\dots,q\}$ and each vector $\nu\in\R_+^d$ with $\nu_1+\dots +\nu_d=1$, we denote by $\Color$ the measure
$P^{\phi,\nu,r}$.

For $\omega\in\{1,\dots,q\}^{\Zd}$, and $\Lambda\subset\Zd$, we note
$n(\Lambda)(\omega)=(n_1(\Lambda)(\omega),\dots ,n_q(\Lambda)(\omega))$, with
$n_k(\Lambda)(\omega)=|\{x\in\Lambda; \omega_x=k\}|$.

\begin{theorem}
\label{empirique}
Let $\phi$ be a measure on $\Sld$ that satisfy to (E),(M) and $\theta_{\phi}=0$
 or $(CLT)$. Let   $q$ be an integer with $q\ge 2$, 
$r\in\{1,\dots,q\}$, and  $\nu\in\R_+^d$ with $\nu_1+\dots +\nu_d=1$.
Then, under $\Color$, we have
$$\frac{n_{\Lambda_t}-|\Lambda_t|((1-\theta_{\phi})\nu+\theta_{\phi}e_{r})}{\sqrt{|\Lambda_t|}}\Longrightarrow\mathcal{N}(0,C),$$
where $C$ is the matrix associated to the quadratic form
$$Q(b)=\chi^f(\phi)(\langle D_\nu b,b\rangle -\langle \nu,b\rangle^2)+\sigma^2_{\phi}\langle e_r-\nu,b\rangle^2,$$
with $D_\nu=\diag(\nu_1,\dots,\nu_d)$.
In other words, $C$ is the matrix of the map
$$b\mapsto \chi^f(\phi)(D_\nu b-\langle \nu,b\rangle \nu)+\sigma^2_{\phi}\langle e_r-\nu,b\rangle (e_r-\nu).$$
\end{theorem}

\begin{proof}
Let $b\in\Rd$ and note 
$Q_t=\frac{n({\Lambda_t})-|\Lambda_t|((1-\theta_{\phi})\nu+\theta_{\phi}e_{r})}{\sqrt{|\Lambda_t|}}$.
For $x\in\Zd$, let us note $Y_x=b_{X_x}$.
We have
\begin{eqnarray*}
\langle n(\Lambda_t),q\rangle & = & \sum_{k=1}^q n_k(\Lambda_t) b_k =  \sum_{k=1}^q \sum_{x\in\Lambda_t}\delta_{X_x}(k) b_k\\
 & = & \sum_{x\in\Lambda_t} \sum_{k=1}^q \delta_{X_x}(k) b_k  =  \sum_{x\in\Lambda_t} b_{X_x}\\
 & = & \sum_{x\in\Lambda_t} Y_x
\end{eqnarray*} 
Now put $m=\langle \nu,b\rangle$ and $z=b_r=\langle e_r,b\rangle$.
We have
 $$\langle Q_n,b\rangle =\frac{\big(\miniop{}{\sum}{x\in\Lambda_t} Y_x\big)-|\Lambda_t|((1-\theta_{\phi})m+\theta_{\phi}z)}{\sqrt{|\Lambda_t|}}.$$
Now if we define $\mu$ to be the image of $\nu$ by $k\mapsto b_k$, it is
not difficult to see that the mean of $\mu$ is $m$ and that 
the law of $(Y_k)_{k\in\Zd}$ under $\Color$ is $P^{\phi,\mu,z}$.
Then, it follows from Theorem~\ref{annealed} that
$\langle Q_n,b\rangle\Longrightarrow\mathcal{N}(0,Q(b))$, with 
$Q(b)=\chi^f(\phi)\sigma^2+(z-m)^2\sigma_{\phi}^2)$, where $\sigma^2$ is the variance of $\nu$. Finally, we get the explicit form
$$Q(b)=\chi^f(\phi)(\langle D_\nu b,b\rangle -\langle \nu,b\rangle^2)+\sigma^2_{\phi}\langle e_r-\nu,b\rangle^2,$$
with $D_\nu=\diag(\nu_1,\dots,\nu_d)$.
Let $L$ be a random vector following $\mathcal{N}(0,C)$, where $C$ is the covariance matrix associated to $Q$.
We have proved that
$$\forall b\in\Rd \langle Q_t,b\rangle\Longrightarrow \langle L,b\rangle.$$
Using the theorem of Levy, it is easy to see that it is equivalent to say 
that $Q_t\Longrightarrow L$.
\end{proof}

We are now interested in having, for $\theta_{\phi}>0$, a version of
theorem~\ref{empirique} in which the observed quantity does not depend on $r$. 
There are several reasons to motivate such a theorem: if we want to use
this central limit theorem to test if a concrete physical 
system conforms to this model (have in mind an Ising or a Potts model for instance), we have a priori no reason to guess the $r$ phase of the underlying theoretical system. There is also a theoretical justification to wish for such a theorem: if we 
get a theorem which ``does not depend'' on $r$, it will be easy to transfer
it to any measure which resides in the convex hull of the  $(\Color)_{1\le r\le q}$.

\begin{theorem}
\label{empirique2}
Let $\phi$ be a measure on $\Sld$ that satisfy to $(E)$, $(M)$ , $\theta_{\phi}>0$ and
$(CLT)$. Let   $q$ be an integer with $q\ge 2$, 
$r\in\{1,\dots,q\}$, and  $\nu\in\R_+^d$ with $\nu_1+\dots +\nu_d=1$.
For $\Lambda\subset\Zd$, we denote
by $R_{\Lambda}$ the element of $\{1,\dots,q\}$ which realizes the maximum
of $(n_{\Lambda}(k)-|\Lambda|(1-\theta_{\phi})\nu(k))_{1\le k\le q}$.
Then, under $\Color$, we have
$$\frac{n_{\Lambda_t}-|\Lambda_t|((1-\theta_{\phi})\nu+\theta_{\phi}e_{R_{\Lambda_t}})}{\sqrt{|\Lambda_t|}}\Longrightarrow\mathcal{N}(0,C),$$
where $C$ is the matrix associated to the quadratic form
$$Q(b)=\chi^f(\phi)(\langle D_\nu b,b\rangle -\langle \nu,b\rangle^2)+\sigma^2_{\phi}\langle e_r-\nu,b\rangle^2,$$
with $D_\nu=\diag(\nu_1,\dots,\nu_d)$.
In other words, $C$ is the matrix of the map
$$b\mapsto \chi^f(\phi)(D_\nu b-\langle \nu,b\rangle \nu)+\sigma^2_{\phi}\langle e_r-\nu,b\rangle (e_r-\nu).$$
\end{theorem}
\begin{proof}
Since $\Color$ is ergodic $\frac{n_{\Lambda_t}}{|\Lambda_t|}=\frac1{|\Lambda_t|}{\sum_{x\in\Lambda_t}}e_{\omega_x}$ almost surely converges to
the mean value of $e_{\omega_0}$, that is $(1-\theta_{\phi})\nu+\theta_{\phi}e_r$. Then, we have the equivalent $n_{\Lambda_t}-|\Lambda_t|(1-\theta_{\phi})\nu\sim |\Lambda_t|\theta_{\phi}e_r$. It follows that $R_{\Lambda_t}=r$ if $t$ is large enough.
Now let $g$ be a bounded continuous function on $\Rd$:
\begin{eqnarray*}
 & & \E g(\frac{n_{\Lambda_t}-|\Lambda_t|((1-\theta_{\phi})\nu+\theta_{\phi}e_{R_{\Lambda_t}})}{\sqrt{|\Lambda_t|}})\\ &  = & \E g(\frac{n_{\Lambda_t}-|\Lambda_t|((1-\theta_{\phi})\nu+\theta_{\phi}e_{r})}{\sqrt{|\Lambda_t|}})\\
& & +\E \big( g(\frac{n_{\Lambda_t}-|\Lambda_t|((1-\theta_{\phi})\nu+\theta_{\phi}e_{R_{\Lambda_t}})}{\sqrt{|\Lambda_t|}})-g(\frac{n_{\Lambda_t}-|\Lambda_t|((1-\theta_{\phi})\nu+\theta_{\phi}e_r)}{\sqrt{|\Lambda_t|}})\big).
\end{eqnarray*}
The first term of the sum converges to the integral of $g$ under $\mathcal{N}(0,C)$ by Theorem~\ref{empirique} and the second one converges to 0 by dominated convergence. It follows that $\E g(\frac{n_{\Lambda_t}-|\Lambda_t|((1-\theta_{\phi})\nu+\theta_{\phi}e_{R_{\Lambda_t}})}{\sqrt{|\Lambda_t|}})$ converges to
the integral of $g$ under $\mathcal{N}(0,C)$ for any bounded continuous function $g$, which is exactly the weak convergence to $\mathcal{N}(0,C)$.
\end{proof}

\begin{theorem}
\label{empirique3}
Let $\phi$ be a measure on $\Sld$ that satisfy to (E),(M),$\theta_{\phi}>0$ and
$(CLT)$. Let   $q$ be an integer with $q\ge 2$, 
$r\in\{1,\dots,q\}$, and  $\nu\in\R_+^d$ with $\nu_1+\dots +\nu_d=1$.
For $\Lambda\subset\Zd$, we denote
by $R_{\Lambda}$ the element of $\{1,\dots,q\}$ which realizes the maximum
of $(n_{\Lambda}(k)-|\Lambda|(1-\theta_{\phi})\nu(k))_{1\le k\le q}$.

Let $\gamma$ be a measure on $\{1,\dots,q\}$ and $\Phi_{\gamma}=\int \Color\ d\gamma(r)$.
Then, under $\Phi_{\gamma}$, we have
$$\frac{n_{\Lambda_t}-|\Lambda_t|((1-\theta_{\phi})\nu+\theta_{\phi}e_{R_{\Lambda_t}})}{\sqrt{|\Lambda_t|}}\Longrightarrow\mu,$$
where $\mu$ is the law of $X+S(e_Z-\nu)$, 
where $X,S$ and $Z$ are independent, with
$X\sim\mathcal{N}(0,C'),S\sim\mathcal{N}(0,\sigma_{\phi}^2)$ and $Z\sim\gamma$.
$C'$  is the matrix associated to the quadratic form
$$Q(b)=\chi^f(\phi)(\langle D_\nu b,b\rangle -\langle \nu,b\rangle^2)$$
with $D_\nu=\diag(\nu_1,\dots,\nu_d)$.
In other words, $C'$ is the matrix of the map
$$b\mapsto \chi^f(\phi)(D_\nu b-\langle \nu,b\rangle \nu).$$
\end{theorem}
\begin{proof}
It just follows from Levy's theorem and a straightforward computation of characteristic function.
\end{proof}

\section{Applications to Potts and Ising models}
In the following, we will always take  
$\phi=\phi^b_{p,q}$, with $b\in\{0,1\}$. We associate to $p$ the inverse 
temperature $\beta=-\frac12\ln(1-p)$. Let us see when the assumptions of
theorem~\ref{empirique} are satisfied.

\begin{itemize}
\item
It is a well-known result that the assumption of ergodicity $(E)$ is always satisfied. 
\item By the inequality of Grimmett and Piza, $(M)$ is always satisfied
when $p<p_g$, or equivalently $\beta<\beta_g=-\frac12\ln(1-p_g)$. (high temperature regime)
\item By the inequality of Pisztora and our lemma~\ref{dimdeux}, $(M)$ holds when $p$ is close enough to 1, or, equivalently, when $\beta$ is large enough. (very low temperature regime)
\item By our theorem~\ref{CLT}, $(CLT)$ holds when $p$ is close enough to 1 ($p>r(q)$, or, equivalently, when $\beta$ is large enough $\beta>\beta_r=-\frac12\ln(1-r(q))$. (very low temperature regime)
\end{itemize}
Remember that when $d=2$, we have proved that $(M)$ and $(CLT)$ hold as soon as $p>r(q)$, which is known to be equal to $p_c(q)$ as soon as $p_g=p_c$.
So, if the conjecture $p_g=p_c$ were proved, we would have the central limit theorem for each value of $\beta$, except for $\beta_c=-\frac12\ln(1-p_c)=\frac12\ln (1+\sqrt q)$.\\

\noindent\textbf{Empirical distributions of Potts models}\\

We have already noticed that if we take $\phi=\phi^1_{p,q}$, with $q\ge 2$ and 
 $\nu=\frac1q(\delta_1+\delta_2+\dots \delta_q)$, $\Color$ is the Gibbs 
measure $\WPt_{q,\beta, r}$ for the $q$-state Potts model on
$\Zd$ at inverse temperature $\beta:= -\frac{1}{2}\log(1-p)$, according
to Proposition~2.4 of Häggström, Jonasson and Lyons~\cite{MR2003f:60173}.
In the high temperature regime $(p<p_c(q))$, we have $\phi^1_{p,q}(0\communique\infty)=0$, so  $\WPt_{q,\beta, r}$ does not depends on $r$.
In fact, it is known that there is uniqueness of the Gibbs measure when
$\phi^1_{p,q}(0\communique\infty)=0$, so $\WPt_{q,\beta, 1}$ is the unique
Gibbs measure at inverse temperature $\beta$. We can now formulate a theorem corresponding to the high temperature regime:
\begin{theorem}
Let $\beta<\beta_g$. There is a unique Gibbs measure for the $q$-state Potts model at inverse temperature $\beta$. If we note $p=1-\exp(-2\beta)$, we have the following results for the empirical distributions: 
$$\frac{n_{\Lambda_t}-|\Lambda_t|\nu}{\sqrt{|\Lambda_t|}}\Longrightarrow\mathcal{N}(0,\frac{\chi(p,q)}{q^2}(qI-J)),$$
where $J$ is the $q\times q$ matrix whose each entry is equal to 1, and
$$\chi(p,q)=\sum_{k=1}^{+\infty} k\phi^1_{p,q}(|C(0)|=k).$$

\end{theorem}

If $p>p_c$, then the Gibbs measures $\phi=\WPt_{q,\beta, r}$ are all distinct (it can be seen as a consequence of the fact that $R_{\Lambda_t}$ 
$\WPt_{q,\beta, r}$-almost surely converges to $r$). Since they are ergodic by theorem~\ref{ergodicite}, they are affinely independent.
Then, in the case $\beta>\beta_r$, we have obtained a central limit theorem for the empirical distribution for a $q$-dimensional convex set of Gibbs measures: the limit is Gaussian when $\phi=\WPt_{q,\beta, r}$, in general it is not Gaussian for a convex combination of them.

\begin{theorem}
\label{basse}
Let $\beta>\beta_r$ and let $\Phi_{\gamma}$ be a Gibbs measure for the $q$ states Potts model at inverse temperature
$\beta$ which can be written in the form  $\Phi_{\gamma}=\int \WPt_{q,\beta, r}\ d\gamma(r)$. 
For $\Lambda\subset\Zd$, we denote
by $R_{\Lambda}$ the element of $\{1,\dots,q\}$ which realizes the maximum
of $(n_{\Lambda}(k))_{1\le k\le q}$. Let us note $p=1-\exp(-2\beta)$.
Then, under $\Phi_{\gamma}$, we have
$$\frac{n_{\Lambda_t}-|\Lambda_t|((1-\theta_{\phi})\nu+\theta_{\phi}e_{R_{\Lambda_t}})}{\sqrt{|\Lambda_t|}}\Longrightarrow\mu,$$
where $\mu$ is the law of $X+S(e_Z-\nu)$, 
when $X,S$ and $Z$ are independent, with
$X\sim\mathcal{N}(0,\frac{\chi^{f}(p,q)}{q^2}(qI-J)),S\sim\mathcal{N}(0,\sigma_{\phi}^2)$ and $Z\sim\gamma$.
 $J$ is the $q\times q$ matrix whose each entry is equal to 1, and
$$\chi^{f}(p,q)=\sum_{k=1}^{+\infty} k\phi^1_{p,q}(|C(0)|=k).$$
\end{theorem}
\noindent{Remark:}
An interesting case of a convex combination is obtain when $\gamma={\nu}$, with the notations of Theorem~\ref{empirique3}. Since we have uniqueness of the infinite cluster in the random cluster cluster, we can consider that
the law of $\Phi_{\nu}$ is obtained by coloring independently each connected
component of the random cluster. Then,  $\Phi_{\nu}$ is just  $\FPt^{\Zd}_{q,\beta}$ in the terminology of Proposition~2.3 of Häggström, Jonasson and Lyons~\cite{MR2003f:60173}.\\

\noindent\textbf{Fluctuations of the magnetization in Ising models}\\

In spite of the fact that $\mu$ is in general not Gaussian, we can observe an intriguing fact when $q=2$, \ie for the Ising model.
In this case $S(e_Z-\nu)=\epsilon S\begin{pmatrix} -\frac12\\ \frac12 \end{pmatrix}$, with $\epsilon=(-1)^{\1_{\{Z=1\}}}S$. But $\epsilon S$ has the same law than $S$. It follows that $\mu$ does not depend on $\gamma$ and is always Gaussian.

Note also that it is known that we have an exponential decay of the covariance
in the Ising model at high temperature -- the exact Ornstein-Zernike directional speed of decay has even be proved by Campanino, Ioffe and Velenik~\cite{MR1964456}. It follows that we have $p_c=p_g$ or equivalently $\beta_c=\beta_g$.
Since the value of the critical point when $d=2$ is the fixed point of 
$x\mapsto x^{*}$, we have even $\beta_r=\beta_c$ when $d=2$.

In this model, it is most relevant to formulate the result 
in term of the magnetization $m_{\Lambda}=n_{\Lambda}.(1 -1)$ rather than
in terms of $n_{\Lambda}$.

\begin{theorem}
\label{basse-ising}
Let $\beta>\beta_r$ and let $\Phi_{\gamma}$ be a Gibbs measure for the Ising model on $\{-1,+1\}^{\Zd}$ at inverse temperature
$\beta$ which can be written in the form  
$\Phi_{\gamma}=\gamma \WPt_{2,\beta, 1}
+(1-\gamma) \WPt_{2,\beta, -1}$.  Let us note $p=1-\exp(-2\beta)$ and
$m_{\Lambda}=\frac1{|\Lambda|}\miniop{}{\sum}{x\in\Lambda}\omega_x$.
Then, under $\Phi_{\gamma}$, we have
$${\sqrt{|\Lambda_t|}}\big(m_{\Lambda_t}-\sign(m_{\Lambda})\theta(p,2)\big)\Longrightarrow\mathcal{N}(0,\chi^f(p,2)+\sigma^2_{p,2}),$$
where
$$\chi^{f}(p,q)=\sum_{k=1}^{+\infty} k\phi^1_{p,q}(|C(0)|=k).$$
Note that $\beta_r=\beta_c$ when $d=2$.
\end{theorem}

Note that for $d=2$, Theorem~\ref{basse-ising} covers the whole set of
Gibbs measure at temperature $\beta>\beta_c$. 
Indeed, $\WPt_{q,\beta, 1}$ and $\WPt_{q,\beta, -1}$
are  known to be the only two extremal Gibbs measures 
when $d=2$. (This celebrated result is due to Higuchi~\cite{MR84m:82020} and Aizenmann~\cite{MR82c:82003}. See also Georgii and Higuchi~\cite{MR2001c:82018} for a modern proof.) 
It follows that every Gibbs measure is a convex combination of
$\WPt_{q,\beta, 1}$ and $\WPt_{q,\beta, -1}$.
We also note that $\theta(p,2)$ appears as the spontaneous magnetization in the ``+'' phase of the Ising model. Since the explicit expression of the spontaneous magnetization is known when $d=2$ -- see Abraham and Martin-L{\"o}f~\cite{MR49:6851}, Aizenman~\cite{MR82c:82003}, and also the bibliographical notes by Georgii~\cite{MR89k:82010} for the whole long story of this result -- , we get for $d=2$ and $p\ge p_c$  the formula   
$\theta(p,2)=\big(1-(\sinh 2\beta)^{-4}\big)^{1/8}=\big(1-16\frac{(1-p)^4}{p^4(2-p)^4}\big)^{1/8}$.

Of course we have a similar theorem in the high temperature regime 
$\beta<\beta_g=\beta_c$.

\begin{theorem}
\label{haute-ising}
Let $\beta<\beta_c$ and let $\Phi$ be a the unique Gibbs measure for the Ising model on $\{-1,+1\}^{\Zd}$ at inverse temperature
$\beta$. We  note $p=1-\exp(-2\beta)$
 and
$m_{\Lambda}=\frac1{|\Lambda|}\miniop{}{\sum}{x\in\Lambda}\omega_x$.

Then, under $\Phi$, we have
$${\sqrt{|\Lambda_t|}}m_{\Lambda_t}\Longrightarrow\mathcal{N}(0,\chi^f(p,2)),$$
where
$$\chi^{f}(p,q)=\sum_{k=1}^{+\infty} k\phi^1_{p,q}(|C(0)|=k).$$
\end{theorem}

Note that for the Gibbs measures $\WPt_{2,\beta, 1}$ or  $\WPt_{2,\beta, -1}$, the central limit theorems could be proved
without the machinery of the above section: since the Ising model satisfy
to the F.K.G. inequalities, it follows from the theorem of Newman that is sufficient to prove that $$\sum_{k\in\Zd}\Cov(\sigma_0,\sigma_k)<+\infty.$$ But it is not difficult
to see that under $\WPt_{2,\beta, 1}$ or  $\WPt_{2,\beta, -1}$,
we have $\Cov(\sigma_0,\sigma_k)=\phi_{p,2}(0\communique k)-\phi_{p,2}(0\communique\infty)^2$
But $\phi_{p,2}(0\communique k)=\phi_{p,2}(0\communique k\text{ by a finite cluster})+\phi_{p,2}(0\communique\infty,k\communique\infty)$.
Then
\begin{eqnarray*}\sum_{k\in\Zd}\Cov(\sigma_0,\sigma_k) & = & \sum_{k\in\Zd}\phi_{p,2}(0\communique k\text{ by a finite cluster})+\sum_{k\in\Zd}\phi_{p,2}(0\communique\infty,k\communique\infty)-\phi_{p,2}(0\communique\infty)^2\\
& = &  \chi^f(p,2)+\sigma^2_{p,2}
\end{eqnarray*}
Nevertheless, when  $\beta>\beta_r$, the Gibbs measure ``with free boundary conditions''  $\FPt^{\Zd}_{2,\beta}$ -- which satisfy to the assumptions of theorem~\ref{basse-ising} -- does not have finite susceptibility: in
this case $$\Cov(\sigma_0,\sigma_k)=\E\sigma_0\sigma_k=\phi_{p,2}(0\communique k)\ge \phi_{p,2}(0\communique\infty,k\communique\infty) \ge\phi_{p,2}(0\communique\infty)^2>0,$$ so the series $\sum_{k\in\Zd}\Cov(\sigma_0,\sigma_k)$ diverges.

These results can be compared with a result of Martin-L{\"o}f~\cite{MR49:8566}: he also proved a central limit theorem for the magnetization  in Ising Models at very low temperature. Particularly, he relays the variance of the limiting normal measure to the second derivative at $0$ of the thermodynamical function $F$.
Nevertheless, his result is slightly different from ours, since he considers
Gibbs measures in large boxes with boundary condition ``+'', whereas we consider here infinite  Gibbs measures under the ``+" phase.

\bibliographystyle{plain}
\bibliography{clust}

\end{document}